\magnification=1200
\vsize=9truein
\hsize=6.5truein
\baselineskip=16pt


  \def\BB{{\cal B}}     
\def\EE{{\cal E}}  \def\FF{{\cal F}}  \def\GG{{\cal G}}

\def\IR{{{\rm I}\!{\rm R}}}

  \def\rightheadline{\tenrm\hfil\folio}
  \def\leftheadline{\tenrm\folio\hfil}

\def\sqr#1#2{{\vcenter{\vbox{\hrule height.#2pt
   \hbox{\vrule width.#2pt height#1pt \kern#1pt
   \vrule width.#2pt}
   \hrule height.#2pt}}}}

\def\rect#1#2#3{\raise .1ex\vbox{\hrule height.#3pt
   \hbox{\vrule width.#3pt height#2pt \kern#1pt\vrule width.#3pt}
        \hrule height.#3pt}}

\def\qed{$\hskip 5pt\rect364$} 

\def\nbar #1{{\raise7pt\hbox{\vrule height2.5pt depth-2pt width6pt}
   \kern-8pt #1}}

\def\nnbar #1{{#1\kern-7.5pt\raise6.2pt\hbox{\vrule height2.3pt depth-2pt
      width6pt}\,}}



\def\ref#1{{\rm [}{\bf #1}{\rm ]}}   
\def\nref#1#2{{\rm [}{\bf #1}{\rm ;\ #2]}}

\def\comp{{\leavevmode
     \raise.2ex\hbox{${\scriptstyle\mathchar"020E}$}}}

\def\today{\ifcase\month\or 
  January\or February\or March\or April\or  
  May\or June\or July\or August\or  
  September\or October\or November\or
  December\fi\space\number\day,\ \number\year}



\outer\def\proclaim#1{\medbreak\noindent\bf\ignorespaces
   #1\unskip.\enspace\sl\ignorespaces}
\outer\def\endproclaim{\par\ifdim\lastskip<\medskipamount\removelastskip
   \penalty 55 \fi\medskip\rm}

\font\tenBbb=msbm10

\newfam\Bbbfam \textfont\Bbbfam=\tenBbb
\font\sevenBbb=msbm7
\scriptfont8=\sevenBbb

\def\picture #1 by #2 (#3){
		\vbox to #2{
				\hrule width #1 height 0pt depth 0pt
				\vfill
				\special{picture #3}}}

\def\scaledpicture #1 by #2 (#3 scaled #4){{
		\dimen0=#1\dimen1=#2
		\divide\dimen0 by 1000\multiply\dimen0 by #4
		\divide\dimen1 by 1000\multiply\dimen1 by #4
		\picture\dimen0 by dimen1 (#3 scaled #4)}}

\outer\def\drop#1{
	\font\cap=Times at 46pt
	\setbox0=\hbox spread 2pt
	{\cap#1\hfil}
	\def\dci{\cap#1}
	\def\dc{\noindent
		\lower34.75pt\hbox
		{\hskip-\hangindent\dci}}
	\dc\vskip-34.75pt\noindent
	\hangindent
	=\wd0\hangafter=-3}

\def\pmb#1{\setbox0=\hbox{#1}%
 \kern-.025em\copy0\kern-\wd0
 \kern.05em\copy0\kern-\wd0
 \kern-.025em\raise.0433em\box0 }


\def\boxit#1{\vbox{\hrule\hbox{\vrule\kern1pt
  \vbox{\kern1pt#1\kern1pt}\kern1pt\vrule}\hrule}}

\newdimen\pxl
\pxl=.2409pt
\newbox\dlbox \setbox\dlbox=\vbox{\hrule width13\pxl height30\pxl
   depth -29\pxl \vskip-38\pxl\hrule width13\pxl height0pt depth1\pxl
   \hbox{\vrule height 30\pxl width1\pxl depth 7\pxl\kern6\pxl\vrule height
   30\pxl width1\pxl depth 7\pxl}}
\newbox\drbox \setbox\drbox=
   \vbox{\hrule width13\pxl height30\pxl depth -29\pxl
   \vskip-38\pxl\hrule width13\pxl height0pt depth 1\pxl\hbox{\kern5\pxl\vrule
   height 30\pxl width1\pxl depth 7\pxl\kern6\pxl\vrule height 30\pxl width1\pxl
   depth 7\pxl}}
\newbox\dlsbox \setbox\dlsbox=\vbox{\hrule width10\pxl height21\pxl depth
   -20\pxl \vskip-27\pxl\hrule width10\pxl height0pt depth 1\pxl\hbox{\vrule
   height 21\pxl width1\pxl depth 5\pxl\kern4\pxl\vrule height 21\pxl
   width1\pxl depth 5\pxl}}
\newbox\drsbox \setbox\drsbox=
   \vbox{\hrule width10\pxl height21\pxl depth -20\pxl
   \vskip-27\pxl\hrule width10\pxl height0pt depth
   1\pxl\hbox{\kern4\pxl\vrule height 21\pxl width1\pxl depth
   5\pxl\kern4\pxl\vrule height 21\pxl width1\pxl depth 5\pxl}}

\def\dl{\,\copy\dlbox\,} 
\def\dr{\,\copy\drbox\,} 
\def\dls{\,\copy\dlsbox\,} 
\def\drs{\,\copy\drsbox\,}
 
\def\dlb{\mathchoice{\dl}{\dl}{\dls}{\dls}}
\def\drb{\mathchoice{\dr}{\dr}{\drs}{\drs}}

\def\si#1#2 #3#4{\if#1[\dlb\else\if#1]\drb\else#1\!#1\fi\fi
   #2,#3\if#4[\dlb\else\if#4]\drb\else#4\!#4\fi\fi}
\def\grph#1#2#3{\if#1[\dlb\else\if#1]\drb\else#1\!#1\fi\fi
   #2\if#3[\dlb\else\if#3]\drb\else#3\!#3\fi\fi}

\def\pagenoslikebook{\nopagenumbers
 \headline={\ifodd\pageno\rightheadline \else\leftheadline\fi}
 \def\rightheadline{\tenrm\hfil\folio}
 \def\leftheadline{\tenrm\folio\hfil}} 

\def\m@th{\mathsurround=0pt }
\def\ialign{\everycr={}\tabskip=0pt \halign}
\def\cases#1{\left\{\,\vcenter{\normalbaselines\m@th 
\ialign{$##\hfil$&\quad##\hfil\crcr#1\crcr}}\right.}

\def\matrix#1{\null\,\vcenter{\normalbaselines\m@th
    \ialign{\hfil$##$\hfil&&\quad\hfil$##$\hfil\crcr
     \mathstrut\crcr\noalign{\kern-\baselineskip}
      #1\crcr\mathstrut\crcr\noalign{\kern-\baselineskip}}}\,}

\def\nin{\noindent}
\def\(#1){{\rm(}#1\/{\rm)}}

\def\nbar #1{{\raise7pt\hbox{\vrule height2.25pt depth-2pt width6pt}
   \kern-8pt #1}}
\def\nbara #1{{\raise4pt\hbox{\vrule height2.25pt depth-2pt width6pt}
   \kern-8pt #1}}

\def\prf{\nin{\it Proof. }}

\def\<{\langle}
\def\>{\rangle}

\def\ov{\overline}\def\reg{\mathop{\rm reg}}

\def\ee{{\bf e}}

\def\BG{1}
\def\CHU{2}
\def\FIT{3}
\def\GET{4}
\def\GS{5}
\def\IMH{6}
\def\ITO{7}
\def\IMK{8}
\def\LOU{9}
\def\MAI{10}
\def\MILa{11}
\def\MILb{12}
\def\PIT{13}
\def\ROG{14}
\def\SAL{15}
\def\SHP{16}
\def\VER{17}
\def\WIL{18}

\centerline{\bf Excursions Above the Minimum for Diffusions
\footnote{*}{\rm This is a
lightly edited form of a manuscript written in the spring of 1985.}}
\medskip
\centerline{P.J. Fitzsimmons}
\centerline{August 23, 2013} 
\bigskip

\nin{\bf 1. Introduction}
\medskip

\nin Let $X=(X_t)_{t\ge 0}$ be a regular diffusion process on an interval
$E\subset\IR$. Let $H_t:=\min_{0\le u\le t}X_u$ denote the past minimum process
of $X$ and consider the excursions of $X$ above its past minimum level:
If $[a,b]$ is a maximal interval of constancy of $t\mapsto H_t$, then $(X_t:a\le
t\le b)$ is the ``excursion above the minimum'' starting at time $a$ and level
$y=H_a$. These excursions, when indexed by the level at which they begin, can be
regarded (collectively) as a point process. The independent increments property
of the first-passage process of $X$ implies that this point process is
Poisonnian in nature, albeit non-homogeneous in intensity. Moreover, intuition
tells us that the distribution of an excursion above the minimum $(X_t:a\le t\le
b)$ should be governed by the It\^o excursion law corresponding to excursions
above the fixed level $y=H_a(\omega)$.

Our first task is to render precise the ruminations of the preceding paragraph.
This is accomplished in sections 2 and 3 by applying Maisonneuve's theory of
exit systems \ref{\MAI} to a suitable auxiliary process $(\overline X_t)$
associated with $X$. The basic result, stated in section 2, affirms the
existence of a ``L\'evy system'' for the point process of excursions of $X$
above its past minimum.

In sections 4, 5, and 6 we discuss several applications of the L\'evy system
constructed in section 3; these applications concern path decompositions of $X$
involving the minimum process $H$. Such decompositions, and related results,
have been found by various authors (see
\ref{\IMH, \LOU, \MILa, \MILb, \ROG, \SAL, \SHP, \VER, \WIL}), most often in the
special case where
$X$ is Brownian motion. The possibility of using L\'evy systems to give a unified
treatment of path decompositions is, of course, not surprising. In an excellent
synthesis \ref{\PIT} Pitman has shown how the existence of a L\'evy system for a
point process attached to a Markov process leads naturally to various path
decompositions of the Markov process.

In section 4 we obtain a general version of Williams' decomposition of a
diffusion at its global minimum. A ``local'' version of Williams' decomposition
can be found in section 5. In section 6 we give a new proof of a result of
Vervaat \ref{\VER}, which states that a Brownian bridge, when split at its
minimum and suitably `` rearranged'' becomes a (scaled) Brownian excursion.
Indeed, we produce an inversion of Vervaat's transformation, showing how a
Brownian excursion may be split and rearranged to yield Brownian bridge.
\bigskip

\nin{\bf 2. Notation and the basic result}
\medskip

\nin Let $X=(\Omega,\FF,\FF_t,\theta_t,X_t,P^x)$ be a canonically defined
regular diffusion on an interval $E\subset\IR$. Here $\Omega$ denotes the space
of paths $\omega\colon[0,+\infty[\to E\cup\{\Delta\}$ which are absorbed in the
cemetery point $\Delta\notin E$ at time $\zeta(\omega)$, and which are
continuous on $[0,\zeta(\omega)[$. For $t\ge 0$, $X_t(\omega)=\omega(t)$, and
$\theta_t\omega$ denotes the path $u\mapsto\omega(u+t)$. The $\sigma$-fields
$\FF$ and $\FF_t$ ($t\ge 0$) are the usual Markovian completions of
$\FF^\circ=\sigma\{X_u:u\ge 0\}$ and $\FF^\circ_t=\sigma\{X_u:0\le u\le t\}$
respectively. The law $P^x$ on $(\Omega,\FF^\circ)$ corresponds to $X$ started
at $x\in E$. We shall also make use of the killing operators $(k_t)$ defined
for $t\ge 0$ by
$$
k_t\omega(u)=\cases{\omega(u),&$u<t$,\cr \Delta,&$u\ge t$.\cr}
$$

Let $A=\inf E$, $B=\sup E$, and write $E^\circ=]A,B[$. We assume throughout the
paper that $A\notin E$, and that $B\in E$ if and only if $B$ is a regular
boundary point which is not a trap for $X$. In particular, these assumptions
imply that the transition kernels of $X$ are absolutely continuous with respect
to the speed measure $m$ (recalled below). See \S 4.11 of It\^o-McKean
\ref{\IMK}.

Let $s$ (resp.\ $m$, resp.\ $k$) denote a scale function (resp.\ speed measure,
resp.\ killing measure) for $X$. Recall from \ref{\IMK} that the generator
$\GG$ of $X$ has the form
$$
\GG u(x)\cdot m(dx)=du^+(x)-u(x)\cdot k(dx),\qquad x\in E^\circ,
\leqno(2.1)
$$
for $u\in D(\GG)$, the domain of $\GG$. Here and elsewhere $u^+$ denotes the
scale derivative:
$$
u^+(x)=\lim_{y\downarrow x}{u(y)-u(x)\over s(y)-s(x)}.
$$

Let $(U^\alpha:\alpha>0)$ denote the resolvent family of $X$. Subsequent
calculations require an explicit expression for the density of $U^\alpha(x,dy)$
with respect to $m(dy)$. Recall from \ref{\IMK} that for each $\alpha>0$ there
are strictly positive, linearly independent solutions $g^\alpha_1$ and
$g^\alpha_2$ of
$$
\GG g(x)=\alpha g(x),\qquad x\in E^\circ;
\leqno(2.2)
$$
$g^\alpha_1$ (resp.\ $g^\alpha_2$) is an increasing (resp.\ decreasing)
solution of (2.2) which also satisfies the appropriate boundary condition at
$A$ (resp.\ $B$). Both $g^\alpha_1$ and $g^\alpha_2$ are uniquely determined up
to a positive multiple. We sometimes drop the superscript $\alpha$, writing
simply $g_1$ and $g_2$. Since $g_1$ and $g_2$ are linearly independent
solutions of (2.2), the Wronskian $W=g_1^+g_2-g_2^+g_1$ is constant. The
resolvent $U^\alpha$ is given by
$$
U^\alpha f(x)=U^\alpha(x,f)=\int_E u^\alpha(x,y)f(y)\,m(dy),
\leqno(2.3)
$$
where
$$
u^\alpha(x,y)=u^\alpha(y,x)=g^\alpha_1(x)g^\alpha_2(y)/W,\qquad x\le y.
\leqno(2.4)
$$
See \S 4.11 of \ref{\IMK} and note that in (2.3) the mass $m(\{B\})$ is the
``stickiness'' coefficient occurring in the boundary condition at $B$ for
elements of $D(\GG)$.

A jointly continuous version $(L^y_t:t\ge 0, y\in E)$ of {\it local time} for
$X$ may be chosen, and normalized to be occupation density relative to $m$, so
that
$$
P^x\int_0^\infty e^{-\alpha t}\,dL^y_t=u^\alpha(x,y).
\leqno(2.5)
$$
Fixing a level $y\in E$, the local time $(L^y_t:t\ge 0)$ is related to the
It\^o excursion law \ref{\ITO}, for excursions from level $y$, as follows. Let
$G(y)$ denote the (random) set of left-hand endpoints (in $]0,\zeta[$) of
intervals contiguous to the level set $\{t>0:X_t=y\}$. Define the hitting time
$T_y$ by
$$
T_y=\inf\{t>0:X_t=y\}\qquad (\inf\emptyset=+\infty).
$$
The It\^o excursion law $n_y$ is determined by the identity
$$
P^x\sum_{u\in G(y)}Z_u\,F\comp k_{T_y}\comp\theta_u=P^x\left(\int_0^\infty
Z_u\,dL^y_u\right)\cdot n_y(F),
\leqno(2.6)
$$
where $x\in E$, $F\in p\FF^\circ$, and $Z\ge 0$ is an $(\FF_t)$-optional
process. Under $n_y$ the coordinate process $(X_t:t> 0)$ is strongly
Markovian with semigroup $(Q^y_t)$ given by
$$
Q^y_t(x,f)=P^x(f\comp X_t;t<T_y).
\leqno(2.7)
$$
The entrance law $n_y(X_t\in\,dz)$ is determined by the corresponding Laplace
transform
$$
W^\alpha f(y)=W^\alpha(y,f)=n_y\int_0^\zeta e^{-\alpha t} f\comp X_t\,dt.
\leqno(2.8)
$$
Conversely, $n_y$ is the unique $\sigma$-finite measure on $(\Omega,\FF^\circ)$
which is carried by $\{\zeta>0\}$ and under which $(X_t:t>0)$ is Markovian with
semigroup (2.7) and entrance law (2.8).

Let $V^\alpha_y$ denote the resolvent of the semigroup $(Q^y_t)$. Taking
$Z_u=e^{-\alpha u}$, $F=\int_0^\zeta e^{-\alpha t}f(X_t)\,dt$ in (2.6), and
using (2.5), we obtain the important identity
$$
U^\alpha f(x)=V^\alpha_yf(x)+u^\alpha(x,y)[m(\{y\})f(y)+W^\alpha f(y)].
\leqno(2.9)
$$
We also recall from \S 4.6 of \ref{\IMK} that the distribution of $T_y$ is
given by
$$
P^x(e^{-\alpha T_y})=\cases{
g^\alpha_1(x)/g^\alpha_1(y),&$x\le y$,\cr g_2^\alpha(x)/g^\alpha_2(y),&$x\ge
y$.\cr}
\leqno(2.10)
$$

Finally, the point process of excursions above the minimum is defined as
follows. For $t\ge 0$ set
$$
\eqalign{
H_t(\omega)&=\cases{\min_{0\le u\le t}X_u(\omega)&if $t<\zeta(\omega)$,\cr
-\infty&if $t\ge\zeta(\omega)$;\cr}\cr
M(\omega)&=\{u>0:X_u(\omega)=H_u(\omega)\};\cr
R_t(\omega)&=\inf\{u>0:u+t\in M(\omega)\};\cr
G(\omega)&=\{u>0: u<\zeta(w), R_{u-}(\omega)=0<R_u(\omega)\}.\cr}
$$
\medskip
\noindent 
Thus $G$ is the random set of left-hand endpoints of intervals
contiguous to the random set $M$. For $u\in G$ we have the excursion ${\ee}^u$
defined by
$$
{\ee}^u_t=\cases{ X_{u+t},&$0\le t<R_u$,\cr \Delta,&$t\ge R_u$.\cr}
$$
The point process $\Pi=({\ee}^u:u\in G)$ admits a L\'evy system as follows. Define
a continuous increasing adapted process $C=(C_t:t\ge 0)$ by
$$
C_t=\cases{s(H_0)-s(H_t),&if $t<\zeta$,\cr C_{\zeta-}&if $t\ge\zeta$.\cr}
$$

\proclaim{(2.11) Theorem} For $Z\ge 0$ and $(\FF_t)$-optional, and $F\in
p\FF^\circ$,
$$
\eqalign{
P^x\sum_{u\in G} Z_u\, F({\ee}^u)&=P^x\int_0^\infty
Z_u\,n^\uparrow_{X_u}(F)\,dC_u\cr
&=P^x\int_A^x Z_{T_y} 1_{\{T_y<+\infty\}} n_y^\uparrow(F)\,ds(y),\cr}
\leqno(2.12)
$$
where $n^\uparrow_y$ denotes the restriction of $n_y$ to
$\{\omega:\omega(t)>y,\forall t\in]0,\zeta(\omega)[\}$.
\endproclaim\bigskip

\nin{\bf (2.13) Remark.} The second equality in (2.12) follows from the first
by the change of variable $u=T_y$. The equality of the first and third terms in
(2.12) amounts to the statement that the time-changed point process
$({\ee}^{T_y}:R_{T_y-}<R_{T_y},A<y<x)$ is a stopped Poisson point process under
$P^x$, with (non-homogeneous) intensity $ds(y)n^\uparrow_y(d\omega)$, stopped at
the first level $y$ for which $T_y=+\infty$. See \ref{\ROG} for this result in
the case of Brownian motion, with or without drift. The general result (2.11)
was suggested by \S 4.10 of \ref{\IMK}.

\bigskip

\nin{\bf 3. Proof of Theorem (2.11)}
\medskip

\nin Maisonneuve's theory of exit systems \ref{\MAI} provides a L\'evy system
description of the point process of excursions induced by a closed, optional,
homogeneous random set. Unfortunately the set $M$ introduced in \S 2 is not
$(\theta_t)$-homogeneous; however the theory of \ref{\MAI} can be brought to
bear once we note that $M$ is homogeneous as a functional of the strong Markov
process $(X_t,H_t)$, $t\ge 0$. This key observation is due to Millar \ref{\MILb}
and has been formalized by Getoor in \ref{\GET}. In the terminology of
\ref{\GET}, the process $H$ is a ``min-functional'': $H_{t+u}=H_t\wedge
H_u\comp\theta_t$. This property ensures that $\ov X:=(X,H)$ is Markovian, as a
simple computation shows.

Following \ref{\GET} we first construct a convenient realization of $\ov X$.
Let $\ov\Omega=\Omega\times(E\cup\{-\infty\})$, $\ov E=\{(x,a)\in E\times
E:a\le x\}$, and for $(\omega,a)\in\ov\Omega$ set
$$
\eqalign{
\ov X_t(\omega,a)&=(X_t(\omega),a\wedge H_t(\omega)),\cr
\ov\theta_t(\omega,a)&=(\theta_t(\omega),a\wedge H_t(\omega)).\cr}
$$
Clearly $\ov X_t\comp\ov\theta_u=\ov X_{t+u}$,
$\ov\theta_t\comp\ov\theta_u=\ov\theta_{t+u}$. Moreover, $M$ can be realized
over $\ov X$ as
$$
\ov M(\omega,a)=\{t>0:\ov X_t(\omega,a)\in D\},
\leqno(3.1)
$$
where $D=\{(x,x):x\in E\}$. Let $\ov\FF^\circ=\sigma\{\ov X_u:u\ge 0\}$, 
$\ov\FF^\circ_t=\sigma\{\ov X_u:0\le u\le t\}$, and for $(x,a)\in\ov E$ let
$\ov P^{x,a}=P^x\otimes\epsilon_a$. The usual Markovian completion of the
filtration $(\ov\FF^\circ_t)$ relative to the laws $(\ov P^{x,a}:(x,a)\in\ov
E)$ is denoted by $(\ov \FF_t)$. Clearly $\ov P^{x,a}(\ov X_0=(x,a))=1$ so that
$\ov X$ has no branch points. Appealing to \S 2 of \ref{\GET} we have the
following

\proclaim{(3.2) Lemma} (i) $\ov X=(\ov\Omega,\ov\FF,\ov\FF_t,\ov\theta_t,\ov
X_t,\ov P^{x,a})$ is a right-continuous, strong Markov process with state space
$\ov E$ and cemetery $\ov\Delta=(\Delta,-\infty)$. The semigroup of $\ov X$
maps Borel functions to Borel functions, so that $\ov X$ is even a Borel right
process.

(ii) Let $\pi\colon(\omega,a)\to\omega$ denote the projection of $\ov\Omega$
onto $\Omega$. If $Z$ is an $(\FF_t)$-optional process, then $Z\comp\pi$ is
$(\ov\FF_t)$-optional.
\endproclaim

Now $\ov M$  is an $(\ov\FF_t)$-optional, $(\ov\theta_t)$-homogeneous set, and
each section $\ov M(\omega,a)$ is closed in $]0,\ov\zeta(\omega,a)[$. Set $\ov
R=\inf\ov M$, so that $\ov R$ is an exact terminal time of $\ov X$ with
$\reg(\ov R)=\{(x,a)\in\ov E:\ov P^{x,a}(\ov R=0)=1\}=D$. This last fact follows
from the regularity of $X$ and the identity
$$
\ov P^{x,a}(\ov R=T_a\comp\pi)=1,\qquad (x,a)\in\ov E.
$$
Let $\ov G$ denote the set of left-hand endpoints of intervals contiguous to
$\ov M$. The properties of the Maisonneuve exit system $({}^*\ov P^{x,a},\ov
K)$ for $\ov M$ are summarized in the next proposition. In what follows, $\ov
\EE^*$ and $\ov\FF^*$ denote the universal completions of $\ov\EE$ (the Borel
sets in $\ov E$) and $\ov\FF^\circ$ respectively.

\proclaim{(3.3) Proposition} {\rm [Maisonneuve]}  There exists a continuous
additive functional (CAF), $\ov K$, of $\ov X$ with a finite 1-potential, and a
kernel
${}^*\ov P^{x,a}$ from $(\ov E,\ov\EE^*)$ to $(\ov\Omega,\ov\FF^*)$ such that
$$
\ov P^{x,a}\sum_{u\in\ov G}\ov Z_u\ov F_u\comp\ov\theta_u=\ov
P^{x,a}\int_0^\infty\ov Z_u\,{}^*\!{\ov P^{\ov X_u}}(\ov F_u)\,d\ov K_u,
\leqno(3.4)
$$
whenever $\ov Z\ge 0$ is $(\ov\FF_t)$-optional and $(u,\ov\omega)\mapsto\ov
F_u(\ov \omega)$ is a $\BB_{[0,+\infty[}\otimes\ov\FF^*$-measurable, positive
function. The CAF $\ov K$ is carried by $D$. For each $(x,a)\in\ov E$, ${}^*
\ov P^{x,a}$ is a $\sigma$-finite measure on $(\ov\Omega,\ov\FF^*)$ under which
the coordinate process is strongly Markovian with the same transition semigroup
as $\ov X$.\endproclaim

\nin{\bf (3.5) Remarks.} The version of $({}^*\ov P^{x,a},\ov K)$ cited in
(3.3) is a variant of that constructed in \ref{\MAI}; the difference stems from 
the possibility that $\ov P^{x,a}(\ov\zeta<+\infty)$ may be positive. The fact
that
$\ov K$ is continuous (and so carried by $D=\reg(\ov R)$) follows from the
construction in \ref{\MAI}, since $\ov M=\{t>0:\ov X_t\in D\}$ and $D$ is
finely perfect (with respect to $\ov X$). Renormalizing the kernel ${}^*\ov
P^{x,a}$ if necessary, we can and do assume that ${}^*\ov P^{y,y}(1-e^{-\ov
R})=1$ for all $y\in E$.

Our plan is to prove Theorem (2.11) by identifying ${}^*\ov P^{x,a}$ and $\ov K$
explicitly, thereby deducing (2.12) from (3.4). First note that by taking $x=y$
in (2.9) and using (2.4) we have
$$
W^\alpha f(y)=\int_{]A,y[}[g_1^\alpha(z)/g_1^\alpha(y)]f(z)\,m(dz)
+\int_{]y,B]}[g_2^\alpha(z)/g_2^\alpha(y)]f(z)\,m(dz),
\leqno(3.6)
$$
where $y\in E$, $\alpha>0$, and $f\ge 0$ is Borel measurable on $E$.

To identify $\ov K$ we define a second CAF of $\ov X$, $\ov C$, by the formula
$$
\ov C_t(\omega,a)=\cases{
s(a\wedge H_0(\omega))-s(a\wedge H_t(\omega))&if $t<\ov\zeta(\omega,a)$,\cr
\ov C_{\ov\zeta-}(\omega,a)&if $t\ge\ov\zeta(\omega,a)$;\cr}
$$
and notice that $\ov C_t(\omega,X_0(\omega))=C_t(\omega)$. Clearly the
fine support of $\ov C$ is $D$. 

For
$x\in E$ put $\psi(x)=W1_{]x,B]}(x)$.

\proclaim{(3.7) Proposition} The CAFs $\ov K$ and $\int_0^t\psi(X_u)\,d\ov
C_u$ are equivalent.\endproclaim

\prf By \nref{\BG}{IV(2.13)} it suffices to check that the CAFs in question
have the same finite 1-potential (over $\ov X$). An argument of Vervaat
\ref{\VER} shows that $P^x(t\in M)=0$ for all $x\in E$ and all $t>0$.
Consequently, $\ov P^{x,a}(t\in\ov M)=0$ for all $(x,a)\in\ov E$ and all $t>0$.
By Fubini's theorem, $\int_0^\infty e^{-t}1_{\ov M}(t)\,dt=0$ a.s.\ $\ov
P^{x,a}$ for all $(x,a)\in\ov E$. Thus taking $\ov Z_u(\ov \omega)=e^{-u}$,
$\ov F_u(\ov\omega)=1-\exp(-\ov R(\ov\omega)\wedge\ov\zeta(\ov\omega))$ in
(3.4), we may compute
$$
\eqalign{
\ov P^{x,a}\int_0^\infty e^{-u}\,d\ov K_u&=\ov P^{x,a}\sum_{u\in\ov G}
e^{-u}\left(\int_0^{\ov R\wedge\ov\zeta} e^{-t}\right)\comp\ov\theta_u\cr
&=\ov P^{x,a}\int_{\ov R\wedge\ov\zeta}^{\ov\zeta} e^{-u}\,du\cr
&=\ov P^{x,a}(e^{-\ov R\wedge\ov\zeta}-e^{-\ov\zeta})\cr
&=P^{x}(e^{-T_a\wedge\zeta}-e^{-\zeta})\cr
&=P^x\int_0^\zeta e^{-t}\,dt-P^x\int_0^{T_a\wedge\zeta} e^{-t}\,dt\cr
&=U^11(x)-V^1_a1(x)\cr
&={u^1(x,a)\over u^1(a,a)}\,U^11(a),\cr}
\leqno(3.8)
$$
where the last equality follows easily from (2.9). On the other hand, our
hypothesis regarding the boundary $A$ implies that $g^1_1(A+)/g_2^1(A+)=0$ (see
\nref{\IMK}{\S 4.6}). Thus
$$
\eqalign{
\ov P^{x,a}\int_0^\infty e^{-t}\psi(X_t)\,d\ov C_t&=P^x\int_{T_a}^\infty
e^{-t}\psi(X_t)\,dC_t\cr
&=P^x\int_A^a e^{-T_y}\psi(y)\,ds(y)\cr
&=\int_A^a[g_2^1(x)/g_2^1(y)]\psi(y)\,ds(y).\cr}
\leqno(3.9)
$$
In (3.9) we have used the change of variables $t=T_y$ to obtain the second
equality, and (2.1) to obtain the third. Now from the definition of the
Wronskian $W$ we see that $d(g_1/g_2)=W\cdot[g_2]^{-2}ds$. Using this
fact and the expression for $\psi$ provided by (3.6) we may continue the
computation begun in (3.9) with
$$
\eqalign{
&=\int_A^a[g_2(x)/g_2(y)]\int_{]y,B]}[g_2(z)/g_2(y)]\,m(dz)\,ds(y)\cr
&=\int_A^a\int_{]y,B]}[g_2(x)g_2(z)/W]\,m(dz)\,d(g_1/g_2)(y)\cr
&=\int_E\int_A^{z\wedge a} d(g_1/g_2)(y)[g_2(x)g_2(z)/W]\,m(dz)\cr
&=\int_E[g_1(z\wedge a)/g_2(z\wedge a)]\cdot[g_2(x)g_2(z)/W]\,m(dz)\cr
&=[u^1(x,a)/u^1(a,a)]U^11(a).\cr}
\leqno(3.10)
$$
The last equality in (3.10) follows from (2.3) and (2.4). In view of
(3.8)--(3.10), we see that $\ov K$ and $\int_0^t\psi(X_s)\,c\ov C_s$
have the same finite 1-potential and so the proposition is proved.\qed\medskip

For $y\in E$ define a measure $\ov Q^y$ on $(\ov\Omega,\ov\FF^\circ)$ by $\ov
Q^y(F)={}^*\ov P^{y,y}(F\comp\ov k_{\ov R})$, where $\ov k_t$ is the killing
operator on $\ov\Omega$. Since ${}^*\ov P^{y,y}(\ov R\not=T_y\comp\pi)=0$, the
first coordinate of $\ov X$, namely $(X_t:t> 0)$, is Markovian under $\ov Q^y$,
with $(Q^y_t)$ as semigroup. Indeed, we claim that $\psi(y)\pi(\ov
Q^y)=n^\uparrow_y$, at least for $ds$-a.e.\ $y\in E$. To verify this claim it
suffices to compare the associated entrance laws.

\proclaim{(3.11) Lemma} Let $f$ be a bounded positive Borel function on $E$.
Then for $ds$-a.e.\ $y\in E$ we have
$$
\psi(y)\ov Q^y\int_0^{\ov\zeta}e^{-\alpha
t}f(X_t)\,dt=W^\alpha f(y),\qquad\forall \alpha>0.
\leqno(3.12)
$$
\endproclaim

\prf Fix $f$ as in the statement of the lemma and also fix $\alpha>0$. For
$y\in E$ write
$$
\gamma(y)=\ov Q^y\int_0^{\ov\zeta} e^{-\alpha t} f(X_t)\,dt.
$$
As noted in the proof of (3.7), we have $\int_0^\infty 1_{\ov M}(t)\,dt=0$,
$\ov P^{x,a}$-a.s.\ for all $(x,a)\in\ov E$. Thus, for $x\in E$,
$$
\eqalign{
U^\alpha f(x)&=\ov P^{x,x}\sum_{u\in\ov G} e^{-\alpha u}\left(\int_0^{\ov R}
e^{-\alpha t} f(X_t)\,dt\right)\comp\ov\theta_u\cr
&=\ov P^{x,x}\int_0^\infty e^{-\alpha u}\,{}^*\!\ov P^{\ov X_u}\left(\int_0^{\ov
R} e^{-\alpha t} f(X_t)\,dt\right)\, d\ov K_u\cr
&=\ov P^{x,x}\int_0^\infty e^{-\alpha u}\gamma(X_u)\psi(X_u)\,d\ov C_u\cr
&=P^x\int_A^x e^{-\alpha T_y}\gamma(y)\psi(y)\,ds(y)\cr
&=\int_A^x[g_2^\alpha (x)/g_2^\alpha (y)]\gamma (y)\psi(y)\,ds(y).\cr}
\leqno(3.13)
$$
On the other hand, by (2.3) and (2.4), we have
$$
U^\alpha
f(x)=\int_{]A,x]}[g_1^\alpha(x)g^\alpha_2(y)/W]f(y)\,m(dy)
+\int_{]x,B]}[g_1^\alpha(y)g_2^\alpha(x)/W]f(y)\,m(dy).
\leqno(3.14)
$$
If we equate the last line displayed in (3.13) with the right side of (3.14),
divide the resulting identity by $g_2^\alpha(x)$, and then differentiate in
$x$, we obtain
$$
\left[\int_{]x,B]}[g_2^\alpha(y)/g_2^\alpha(x)]f(y)\,m(dy)\right]\,ds(x)
=\gamma(x)\psi(x)\,ds(x)
$$
as measures on $E$, and the lemmas follows.\qed\medskip

\proclaim{(3.15) Corollary} For $ds$-a.e.\ $y\in E$, $\psi(y)\pi(\ov
Q^y)=n^\uparrow_y$ as measures on $(\Omega,\FF^\circ)$.\endproclaim

\prf As noted earlier, both $\psi(y)\pi(\ov Q^y)$ and $n^\uparrow_y$ make the
coordinate process on $(\Omega,\FF^\circ)$ into a Markov process with
transition semigroup $(Q^y_t)$. It follows from Lemma (3.11) that these
measures have the same one-dimensional distributions (and consequently the same
finite dimensional distributions) for $ds$-a.e.\ $y$. Since
$\FF^\circ=\sigma(X_u:u\ge 0)$ is countably generated, the corollary
follows.\qed\medskip

\nin {\it Proof of Theorem (2.11).} Let $Z\ge 0$ be $(\FF_t)$-optional and let
$F\ge 0$ be $\FF^\circ$-measurable. By (3.2)(ii), the process $Z\comp\pi$ is
$(\ov\FF_t)$-optional. We may now use (3.4), (3.7), and (3.15) to compute
$$
\eqalign{
P^x\sum_{u\in G}Z_u F({\ee}^u)
&=\ov P^{x,x}\sum_{u\in\ov G}Z_u\comp\pi F(\pi\comp \ov k_{\ov
R}\comp\ov\theta_u)\cr
&=\ov P^{x,x}\int_0^\infty Z_u\comp\pi\,\,{}^*\ov P^{\ov X_u}(F(\pi\comp\ov k_{\ov
R}))\,d\ov K_u\cr
&=\ov P^{x,x}\int_0^\infty Z_u\comp\pi\,\,\ov Q^{\ov X_u}(F(\pi\comp\ov k_{\ov
R}))\,\psi(X_u)\,d\ov C_u\cr
&= P^x\int_0^\infty Z_u\,n^\uparrow_{X_u}(F)\,dC_u.\cr}
$$
The proof of Theorem (2.11) is  complete.\qed
\bigskip

\nin{\bf 4. Williams' decomposition}
\medskip

\nin In this section we use the L\'evy system (2.12) to obtain a new proof (of
a general version) of Williams' decomposition \ref{\WIL} of a diffusion at its
global minimum. A more ``computational'' proof of Williams' theorem, based on
the same idea used in the present paper, may be found in \ref{\FIT}.

For simplicity we assume that $\gamma:=H_{\zeta-}$ satisfies $P^x(\gamma>A)=1$
for all $x\in E$. We also assume that $\rho:=\inf\{t>0:X_t=\gamma\}$ satisfies
$P^x(\rho<\zeta)=1$ for all $x\in E$. Then $\rho$ is the unique time at which
$X$ takes its global minimum value $\gamma$ ({\it cf.} \ref{\VER}). Note that
for $x\ge y$ (both in $E$),
$$
P^x(\gamma>y)=P^x(T_y=+\infty).
\leqno(4.1)
$$
Define a function $r$ on $E$ by
$$
r(x)=\cases{
P^x(T_{x_0}<+\infty),&$x\ge x_0$,\cr
[P^{x_0}(T_x<+\infty)]^{-1},&$x<x_0$,\cr}
$$
where $x_0\in E^\circ$ is fixed but arbitrary. Clearly $r$ is strictly positive
and decreasing. Arguing as in \nref{\IMK}{pp.~128--129} one may check that $r$ is
the unique positive decreasing solution of $\GG r\equiv 0$ on $E^\circ$ which
satisfies $r(x_0)=1$ and the boundary condition at $B$. Note that
$$
P^x(T_y<+\infty)=r(x)/r(y),\qquad x>y.
\leqno(4.2)
$$

Before proceeding to the decomposition theorem we need a preliminary result.

\proclaim{(4.3) Lemma} For $y\in E^\circ$ let $S_y=\inf\{t>0:X_{t-}=y\}$. Then
$$
n_y^\uparrow (S_y=+\infty)=-{r^+(y)\over r(y)},\qquad\forall y\in E^\circ.
$$
(Recall that $r^+=d^+r/ds^+$.)\endproclaim

\prf Let $q$ be an increasing solution of $\GG q\equiv 0$ on $E^\circ$ such
that $q$ is linearly independent of $r$. We assume that $q$ is normalized so
that the Wronskian $q^+r-r^+q$ is identically 1. Fix $a<b$ both in $E^\circ$
and let $v_{ab}$ denote the potential density (relative to $m$) of $X$ killed
at time $T_a\wedge T_b$. One checks that for $x\le y$,
$$
v_{ab}(x,y)=v_{ab}(y,x)={D(a,x)D(y,b)\over D(a,b)},
$$
where $D(x,y)$ is the determinant
$$
\left|\matrix{q(y)&q(x)\cr r(y)&r(x)\cr}\right|.
$$
Note that $D(x,y)>0$ if $x<y$. Now let $y\in]a,b[$ and use (2.6) to compute
$$
\eqalign{
P^y(T_b<T_a)&=P^y\sum_{u\in G(y)}1_{\{u<T_a\wedge
T_b\}}1_{\{T_b<T_y\}}\comp\theta_u\cr
&=P^y(L^y_{T_a\wedge T_b})\,n_y(T_b<\zeta).\cr}
$$
But clearly $P^y(T_b<T_a)=D(a,y)/D(a,b)$ while $P^y(L^y_{T_a\wedge
T_b})=v_{ab}(y,y)$, so that
$$
n_y(T_b<\zeta)=[D(a,y)/D(a,b)]/v_{ab}(y,y)=D(y,b)^{-1}.
$$
Finally,
$$
\eqalign{
n^\uparrow_y(S_y=+\infty)&=\lim_{x\downarrow y}n^\uparrow_y(T_x<+\infty,
S_y=+\infty)\cr
&=\lim_{x\downarrow y}n^\uparrow_y(T_x<+\infty)\,P^x(T_y=+\infty)\cr
&=\lim_{x\downarrow y}n_y(T_x<\zeta)[1-r(x)/r(y)]\cr
&=\lim_{x\downarrow y}\left[{1\over r(y)}\cdot{r(y)-r(x)\over s(x)-s(y)}\cdot
{s(x)-s(y)\over D(y,x)}\right]\cr
&=-r(y)^{-1}\cdot r^+(y),\cr}
$$
since $\lim_{x\downarrow y}[s(x)-s(y)]/D(y,x)$ is the reciprocal of the
Wronskian $q^+r-r^+q\equiv 1$.\qed\medskip

Now define probability laws on $(\Omega,\FF^\circ)$ by
$$
P^{x\downarrow}_y(F)=P^x(F\comp k_{T_y}|T_y<+\infty),
\leqno(4.4)
$$
$$
P^\uparrow_y(F)=n^\uparrow_y(F|S_y=+\infty),
\leqno(4.5)
$$
whenever $x>y>A$. The coordinate process is a diffusion under any of these
laws: $P^{x\downarrow}_y$ is the law of $X$ started at $x$, conditioned to
converge to $A$, and then killed at $T_y$; $P^\uparrow_y$ is the law of $X$
started at $y$ and conditioned to never return to $y$. These conditionings are
accomplished by means of the appropriate $h$-transforms. In particular, the
associated infinitesimal generators are given by
$$
\GG^{x\downarrow}_yf(z)=r(z)^{-1}\GG(fr)(z),\qquad z>y;
\leqno(4.6)
$$
$$
\GG^\uparrow_yf(z)=r_y(z)^{-1}\GG(fr_y)(z),\qquad z>y,
\leqno(4.7)
$$
where $r_y(z)=1-r(z)/r(y)$.

We can now state the general version of Williams' theorem. Recall that
$\gamma=H_{\zeta-}$ and $\rho=\inf\{t>0:X_t=\gamma\}$.

\proclaim{(4.8) Theorem} (a) The joint law of $(\gamma,\rho,\zeta)$ is given by
$$
P^x(f(\gamma)e^{-\alpha\rho-\beta\zeta})=
\int_A^x[g_2^{\alpha+\beta}(x)/g_2^{\alpha+\beta}(y)]f(y)P^\uparrow_y(e^{-\beta\zeta}){-dr(y)\over
r(y)}.
\leqno(4.9)
$$

(b) For $F,G\in b\FF^\circ$ and $\psi$ bounded and Borel on $E$,
$$
P^x(F\comp
k_\rho\psi(\gamma)G\comp\theta_\rho)=P^x(P^{x\downarrow}_\gamma(F)\psi(\gamma)P^\uparrow_\gamma(G)).
\leqno(4.10)
$$
\endproclaim

\nin{\bf (4.11) Remark.} The intuitive content of (4.10) is that the processes
$(X_t:0\le t<\rho)$ and $(X_{\rho+t}:0\le t<\zeta-\rho)$ are conditionally
independent under $P^x$, given $\gamma$; and that the conditional
distributions, given that $\gamma=y$, are $P^{x\downarrow}_y$ and $P^\uparrow_y$
respectively.
\medskip

\nin{\it Proof of (4.8).} Define $J(y,\omega)=1_{\{S_y=+\infty\}}(\omega)$ and
observe that $\rho(\omega)=u$ if and only if $u\in G(\omega)$ and
$J(X_u(\omega),{\ee}^u(\omega))=1$. Thus, using (2.11),
$$
\eqalign{
P^x(F\comp k_\rho\psi(\gamma)G\comp \theta_\rho)&=
P^x\sum_{u\in G}F\comp k_u\psi(X_u)G({\ee}^u)J(X_u,{\ee}^u)\cr
&=\int_A^x P^x(F\comp T_y;T_y<+\infty)\psi(y) n^\uparrow_y(G\cdot
J(y,\cdot))\,ds(y)\cr
&=\int_A^x P^{x\downarrow}(F)
\psi(y)P^\uparrow_y(G)P^x(T_y<+\infty)n^\uparrow_y(S_y=+\infty)\,ds(y).\cr}
\leqno(4.12)
$$
Taking $F=G=1$ in (4.12) we see that
$$
P^x(\gamma\in dy)=P^x(T_y<+\infty)n^\uparrow_y(S_y=+\infty)\,ds(y).
\leqno(4.13)
$$
Now (4.13) substituted into the last line of (4.12) yields (4.10). To obtain
(4.9) use (4.10) with $F=e^{-(\alpha+\beta)\zeta}$ and $G=e^{-\alpha\zeta}$,
noting that $F\comp k_\rho=e^{-(\alpha+\beta)\rho}$ and
$\zeta=\rho+\zeta\comp\theta_\rho$ ($P^x$-a.s.) since $\rho<\zeta$, $P^x$-a.s.
Thus
$$
\eqalign{
P^x(f(\gamma)e^{-\alpha\rho}e^{-\beta\zeta})
&=P^x(f(\gamma)[e^{-(\alpha+\beta)\zeta}]\comp
k_\rho\, [e^{-\beta\zeta}]\comp\theta_\rho)\cr
&=P^x(P^{x\downarrow}_\gamma(e^{-(\alpha+\beta)\zeta})
f(\gamma)P^\uparrow_\gamma(e^{-\beta\zeta})).\cr}
$$
Formula (4.9) now follows since 
$$
\eqalign{
P^{x\downarrow}_y(e^{-(\alpha+\beta)\zeta})&=
P^x(e^{-(\alpha+\beta)T_y})/P^x(T_y<+\infty)\cr
&=[g_2^{\alpha+\beta}(x)/g_2^{\alpha+\beta}(y)]/P^x(T_y<+\infty),\cr}
$$
and since $n^\uparrow_y(S_y=+\infty)=-r^+(y)/r(y)$ (Lemma (4.3)).\qed\medskip

\proclaim{(4.14) Corollary} $P^x(\rho\in dt, \gamma\in dy)=P^x(T_y\in
dt){\displaystyle{-dr(y)\over r(y)}}$.
\endproclaim
\bigskip

\nin{\bf 5. A local decomposition}
\medskip

\nin Fix $t>0$ and define
$$
\rho_t=\inf\{u>0:X_u=H_t\}\wedge t.
$$
Arguing as in \ref{\VER} one can show that, almost surely on $\{t<\zeta\}$,
$\rho_t$ is the unique $u\in]0,t[$ such that $X_u=H_t$. Our purpose in this
section is to describe the conditional distribution of $\{X_u:0\le u\le t\}$
under $P^b$, given that $H_t=y$, $\rho_t=u$, and $X_t=x$. This conditional
distribution has been computed by Imhof \ref{\IMH} for the Brownian motion (and
closely related processes). The joint law of $(H_t,\rho_t,X_t)$, again in the
case of Brownian motion, has been found by Shepp \ref{\SHP}. See also
\ref{\CHU, \LOU, \SAL} for related results.

We begin by computing the joint law of $(H_t, \rho_t, X_t)$. Recall from
\nref{\IMK}{\S 4.11} that the first passage distribution $P^x(T_y\in dv)$ has a
density $f(v;x,y)$ on $]0,+\infty[$ relative to Lebesgue measure. Note that if
we set $F_{t,y}(x)=P^x(t<T_y<+\infty)$, then (see \nref{\IMK}{p.\ 154})
$$
f(t;x,y)=-{\partial\over\partial t}F_{t,y}(x)=\GG F_{t,y}(x),\qquad x>y\in E, t>0.
\leqno(5.1)
$$
Applying $Q^y(z,dx)$ to both sides of (5.1) and integrating over
$x\in]y,+\infty[\cap E$ (making use of the relation $Q^y_s\GG=\GG Q^y_s$ on
$]y,+\infty[$), we obtain
$$
f(t+s; z,y)=\int_{]y,+\infty[} Q^y_s(z,dx)f(t;x,y).
$$
In other words, $(t,x)\mapsto f(t;x,y)$ is an exit law for the semigroup
$(Q^y_s)$.

Next, recall from \nref{\IMK}{\S 4.11} that the semigroup $(Q^y_t)$ has a
density $q^y(t;x,z)$ (for $x\wedge z>y$) relative to the speed measure $m(dz)$;
we have $q^y>0$ on $]0,+\infty[\times(]y,B[)^2$ and $q^y(t;x,z)=q^y(t;z,x)$. The
entrance law for $n^\uparrow_y$ can now be expressed as
$$
n^\uparrow_y(X_t\in dx)=q^\uparrow_y(t;x)\,m(dx),
\leqno(5.2)
$$
where
$$
q^\uparrow_y(t;x)=\int_{]y,B]} n^\uparrow_y(X_{t-u}\in dz)q^y(u;z,x).
\leqno(5.3)
$$
Substituting (5.2) into (5.3) and using the symmetry of $q^y$, we see that
$$
q^\uparrow_y(t+u;x)=\int_{]y,B]} Q^y_u(x,dz)q^\uparrow_y(t;z).
\leqno(5.4)
$$
But (5.4) means that $(t,x)\mapsto q^\uparrow_y(t;x)$ is also an exit law for
$(Q^y_t)$. Finally, using (3.6), if $\alpha>0$ and $h$ is positive, measurable,
and vanishes off $]y,B]$, we may compute
$$
\eqalign{
\int_0^\infty e^{- \alpha t}\int_Eq^\uparrow_y(t;x) h(x)\,m(dx)\,dt
&=W^\alpha h(y)\cr
&=\int_{]y,B]}[g_2^\alpha(x)/g_2^\alpha(y)]h(x)\,m(dx)\cr
&=\int_{]y,B]} P^x(e^{-\alpha T_y}) h(x)\,m(dx)\cr
&=\int_0^\infty e^{-\alpha t}\int_E f(t;x,y) h(x)\,m(dx).\cr}
$$
By Laplace inversion,
$$
q^\uparrow_y(t;x)=f(t;x,y)
\leqno(5.5)
$$
for $dt\otimes dm$-a.e.\ $(t,x)$ in $]0,+\infty[\times]y,B]$. Since both sides
of (5.5) are exit laws (and so excessive functions in time-space), it follows
that (5.5) holds identically for $t>0$, $y\in E^\circ$ and $E\ni x>y$. See \S 3
of \ref{\GS}, and especially (3.17) therein.

\proclaim{(5.6) Proposition} For $b\in E$, $x\in E$, $u\in ]0,t[$, and $y\in
]A,b\wedge x[$,
$$
P^b(H_t\in dy, \rho_t\in du, X_t\in dx)=f(u;b,y)f(t-u;x,y)\,ds(y)\,du\,m(dx).
\leqno(5.7)
$$
\endproclaim

\prf Let $g$, $h$, and $\phi$ be bounded positive Borel functions on $\IR$,
vanishing off $E$, $E$, and $]0,t[$, respectively. Put
$J(v,y,\omega)=1_{\{S_y>v\}}(\omega)$. Using (2.11) we have, since
$u=\rho_t(\omega)$ if and only if $u\in G(\omega)$ and
$J(t-u,X_u(\omega),{\ee}^u(\omega))=1$,
$$
\eqalign{
P^b(g(H_t)\phi(\rho_t)h(X_t))
&=P^b\sum_{u\in G}g(X_u)\phi(u)h(X_{t-u}\comp\theta_{u})J(t-u, X_u, {\ee}^u)\cr
&=P^b\sum_{u\in G}g(X_u)\phi(u) h({\ee}^u_{t-u})J(t-u, X_u, {\ee}^u)\cr
&=\int_A^b\int_\Omega g(y)\phi(T_y(\omega))
n^\uparrow_y(h(X_{t-T_y(\omega)})P^b(d\omega)\,ds(y)\cr
&=\int_A^b\int_{]0,t[} g(y)\phi(u) n^\uparrow_y(h(X_{t-u})
f(u;b,y)\,du\,ds(y).\cr}
$$
The proposition now follows from (5.2) and (5.5).\qed\medskip

Our local decomposition of $X$ will be expressed in terms of certain ``bridges''
of $X$. First, let $\hat K^{y,\ell,x}$ denote the $h$-transform of
$P^\uparrow_y$ by means of the time-space harmonic function
$$
h_{\ell,x}(t,z)=q^y(\ell-t;z,x)\left[{r(y)-r(x)\over r(y)-r(z)}\right]
1_{]0,\ell[}(t),
$$
where $\ell>0$ and $x>y$. Straightforward computations show that the absolute
probabilities and transition probabilities under $\hat K^{y,\ell,x}$ are given by
$$
\hat K^{y,\ell,x}(X_t\in dz)={q^y(\ell-t;z,x)f(t;y,z)\over f(\ell;y,x)}\,m(dz),
$$
and
$$
\hat K^{y,\ell,x}(X_{t+v}\in dw|X_t=z)={q^y(v;z,w)q^y(\ell-t-v;w,x)\over
q^y(\ell-t;z,x)}\,m(dw).
$$
Moreover ({\it cf.} \ref{\SAL})
$$
\hat K^{y,\ell,x}(\zeta=\ell, X_{\zeta-}=x)=1,
\leqno(5.8)
$$
$$
\int_{]y,B]} \hat K^{y,\ell,x}(F)\, P^\uparrow_y(X_\ell\in
dx)=P^\uparrow_y(F\comp k_\ell).
$$
Thus, $\{\hat K^{y,\ell,x}:x\in]y,B]\}$ is a regular version of the conditional
probabilities $F\mapsto P^\uparrow_y(F\comp k_\ell|X_\ell=x)$.

Now let $K^{x,\ell,y}$ denote the image of $\hat K^{y,\ell,x}$ under the
time-reversal mapping, taking $\omega$ to the path $\gamma_\ell\omega$ defined by
$$
(\gamma_\ell\omega)(t)=\cases{
\omega(\ell-t),&$0<t<\ell$\cr
\omega(\ell-),&$t=0$\cr
\Delta,&$t\ge \ell$.\cr}
$$
Like $\hat K^{y,\ell,x}$, $K^{x,\ell,y}$ is the law of a non-homogeneous Markov
diffusion; from (5.8) we see that
$$
K^{x,\ell,y}(X_0=x,\zeta=\ell,X_{\zeta-}=y)=1.
$$
Moreover, computation of finite dimensional distributions shows that the
transition probabilities for $K^{x,\ell,y}$ are given by
$$
K^{x,\ell,y}(X_{t+v}\in dw|X_t=z)={q^y(v;z,w) f(\ell-t-v;w,y)\over
f(\ell-t;z,y)}.
\leqno(5.9)
$$
It follows that $\{K^{x,\ell,y}:\ell>0\}$ is a regular version of the
conditional probabilities
$$
P^{x\downarrow}_y(\cdot|\zeta=\ell).
$$

\proclaim{(5.10) Theorem} Let $b\in E$. Then under $P^b$ the path fragments
$(X_t: 0\le t<\rho_t)$ and $(X_{\rho_t+u}:0\le u<t-\rho_t)$ are conditionally
independent given $(H_t,\rho_t,X_t)$ on $\{X_t\in E\}=\{t<\zeta\}$. Moreover,
given that $H_t=y$, $\rho_t=u$, and $X_t=x$ ($0<u<t$, $y>x$), the above processes
have conditional laws $K^{b,u,y}$ and $\hat K^{y,t-u,x}$ respectively.
\endproclaim

The proof of (5.10) is similar to that of (4.8) and is left to the interested
reader as an exercise.
\bigskip

\nin{\bf 6. A result of W.\ Vervaat}

\nin In this last section we use the decomposition of \S 5 to give a new proof
of a result of Vervaat \ref{\VER} which concerns a path transformation carrying
Brownian bridge into Brownian excursion.

In this section we take the basic process $(X_t,P^x)$ to be standard Brownian
motion on $\IR$. Let $P_0$ denote the law of {\it Brownian bridge}; namely, 
$$
P_0(F)=P^0(F|X_1=0),\qquad F\in\FF_1.
$$
Under $P_0$ the coordinate process is centered Gaussian with continuous paths,
$X_0=0$, and covariance $P_0(X_uX_t)=u(1-t)$ for $0\le u\le t\le 1$.

Next, Let $P_+$ denote the law of scaled Brownian excursion. Under $P_+$ the
coordinate process $(X_t:0\le t\le 1)$ is a non-homogeneous Markov diffusion
with absolute probabilities 
$$
P_+(X_t\in dx)={2x^2\over \sqrt{2\pi t^3(1-t)^3}}e^{-x^2/2t(1-t)}
\leqno(6.1)(a)
$$
and transition probabilities
$$
P_+(X_{t+v}\in dy|X_t=x)=p(v;y-x)\left({1-t\over
1-t-v}\right)^{3/2}{y\exp(-y^2/2(1-t-v))\over x\exp(-x^2/2(1-t))},
\leqno(6.1)(b)
$$
where $p(v;x)=(2\pi v)^{-1/2}e^{-x^2/2v}$ is the Gauss kernel, and $0<t<t+v<1$,
$0<x,y$. Also, $P_+(\zeta=1)=P_+(X_t>0,\forall t\in]0,1[)=P_+(X_0=X_{1-}=0)=1$.

Computation of finite dimensional distributions now shows that the following
identities hold:
$$
\eqalign{
k_u(P_+(\cdot|x_u=y))&=\hat K^{0,u,y},\cr
\theta_u(P_+(\cdot|X_u=y))&=K^{y,1-u,0},\cr}
$$
where $\hat K^{0,u,y}$ and $K^{y,1-u,0}$ are as defined in the last section, the
basic process being standard Brownian motion.

Now let $\Omega_0=\{\omega\in\Omega:\omega(0)=\omega(1-)=0,
\zeta(\omega)=1\}$ and $\overline \Omega=\Omega_0\times]0,1[$. Define a map
$\Phi\colon
\overline
\Omega\to\Omega_0$ by
$$
\Phi(\omega,u)(t)=\Phi_u(\omega)(t)=\cases{
\omega(u+t)-\omega(u),&$0\le t<1-u$,\cr
\omega(u+t-1)-\omega(u),&$1-u\le t<1$.\cr}
$$
In the following we regard $P_+$ and $P_0$ as measures on $\Omega_0$. Define
$\overline P$ on $\overline\Omega$ by $\overline P=P_+\otimes\lambda$, where
$\lambda$ is Lebesgue measure on $]0,1[$. Set $U(\omega,u)=u$ and $V=1-U$ on
$\overline\Omega$.

\proclaim{(6.2) Proposition} The joint law of $(\Phi,V,X_U)$ under $\overline P$
is the same as the joint law of $(\omega,\rho_1,-H_1)$ under $P_0$.\endproclaim

\prf For paths $\omega$ and $\omega'$, and $t\in]0,1[$ let $\omega/t/\omega'$
denote the spliced path
$$
(\omega/t/\omega')(u)=\cases{
\omega(u),&$0\le u<t$,\cr
\omega'(u-t),&$t\le u<1$,\cr}
$$
and let $\tau_y\omega(t)=\omega(t)-y$. Let $p_+(u,y)=P_+(X_u\in dy)/dy$. Note
that if $\omega(u)=\omega'(0)$, then
$$
\Phi(\omega/u/\omega',u)=(\tau_y\omega'/1-u/\tau_y\omega),
$$
where $0<u<1$ and $y=\omega(u)$. Thus,
$$
\eqalign{
\overline P(&F\comp\Phi\,\psi(V,X_U))
=\int_0^1 P_+(F\comp\Phi_u\,\psi(1-u,X_u))\,du\cr
&=\int_0^1\int_0^\infty\int_\Omega\int_\Omega
  F(\Phi_u(\omega/u/\omega'))\psi(1-u,y)\hat
K^{0,u,y}(d\omega)K^{y,1-u,0}(d\omega') p_+(u,y)\,dy\,du\cr
&=\int_0^1\int_0^\infty\int_\Omega\int_\Omega
  F(\tau_y\omega'/1-u/\tau_y\omega)\psi(1-u,y)\hat
K^{0,u,y}(d\omega)K^{y,1-u,0}(d\omega') p_+(u,y)\,dy\,du\cr
&=\int_0^1\int_0^\infty\int_\Omega\int_\Omega
  F(\omega'/1-u/\omega)\psi(1-u,y) K^{0,1-u,-y}(d\omega') K^{-y,u,0}(d\omega)
p_+(u,y)\,dy\,du\cr
&=P_0(F\cdot\psi(\rho_1,-H_1)).\cr}
$$
\qed\medskip

\proclaim{(6.3) Corollary} (Vervaat): Define a transformation
$\Psi:\Omega_0\to\Omega_0$ by
$$
(\Psi\omega)(t)=\cases{
\omega(\rho_1(\omega)+t)-H_1(\omega),&$0\le t<1-\rho_1(\omega)$,\cr
\omega(\rho_1(\omega)+t+1)-H_1(\omega),&$1-\rho_1(\omega)\le t<1$.\cr}
$$
Then $\Psi(P_0)=P_+$. That is, the $P_0$-law of $(X_t\comp\Psi:0\le t<1)$ is
$P_+$.\endproclaim

\prf It is easy to check that $\Psi\comp\Phi(\omega,u)=\omega$ for all
$(\omega,u)\in\overline\Omega$. Using Proposition (6.2),
$$
\eqalign{
P_0(F\comp\Psi)
&=\overline P(F\comp\Psi\comp\Phi)\cr
&=\overline P(F\comp\pi_1)\cr
&=P_+(F),\cr}
$$
where $\pi_1:(\omega,u)\to\omega$.\qed\bigskip\bigskip

\baselineskip=13pt
\frenchspacing

\centerline{\bf References}
\bigskip

\itemitem{[\BG]}
Blumenthal, R.M. and  Getoor, R.K.:
{\it Markov Processes and Potential Theory.}
Academic Press, New York-London 1968.
\smallskip

\itemitem{[\CHU]}
Chung, K.L.: Excursions in Brownian motion. {\it Ark. Mat.} {\bf 14} (1976)
155--177.
\smallskip

\itemitem{[\FIT]}
Fitzsimmons, P.J.: Another look at Williams' decomposition theorem. {\it Seminar
on Stochastic Processes, 1985\/}, pp~79--85, Birkh\"auser
Boston, 1986. 
\smallskip

\itemitem{[\GET]}
Getoor, R.K.: Splitting times and shift functionals. {\it Z. Wahrsch. verw.
Gebiete} {\bf 47} (1979) 69--81.
\smallskip

\itemitem{[\GS]}
Getoor, R.K. and Sharpe, M.J.: Excursions of dual processes. {\it Adv. Math. }
{\bf 45} (1982) 259--309.
\smallskip

\itemitem{[\IMH]}
Imhof, J.-P.:
Density factorizations for Brownian motion, meander and the three-dimensional Bessel process,
and applications. 
{\it J. Appl. Probab.} {\bf  21} (1984) 500--510. 
\smallskip

\itemitem{[\ITO]}
It\^o, K.: Poisson point processes attached to Markov processes. {\it Proceedings
of the Sixth Berkeley Symposium on Mathematical Statistics and Probability} (Univ.
California, Berkeley, Calif., 1970/1971), Vol. III, pp.
225--239. Univ. California Press, Berkeley, 1972.
\smallskip

\itemitem{[\IMK]}
It\^o, K. and McKean, H.P.:  {\it Diffusion Processes and their Sample Paths.}
(Second printing, corrected.) Springer-Verlag, Berlin-New York, 1974. 
\smallskip

\itemitem{[\LOU]}
Louchard, G.: Kac's formula, L\'evy's local time and Brownian excursion. {\it J.
Appl. Probab.} {\bf 21} (1984)  479--499.
\smallskip

\itemitem{[\MAI]} 
Maisonneuve, B.:
Exit systems.
{\it Ann. Probab.} {\bf 3} (1975) 399--411. 
\smallskip

\itemitem{[\MILa]}
Millar, P.W.: Zero-one laws and the minimum of a Markov process. {\it Trans.
Amer. Math. Soc.\/} {\bf 226} (1977) 365--391.
\smallskip

\itemitem{[\MILb]}
Millar, P.W.: A path decomposition for Markov processes. {\it Ann. Probab.\/}
{\bf  6} (1978)  345--348.
\smallskip

\itemitem{[\PIT]}
 Pitman, J.W.: L\'evy systems and path decompositions. {\it Seminar on
Stochastic Processes, 1981}, pp. 79--110,
Birkh\"auser, Boston, 1981. 
\smallskip

\itemitem{[\ROG]}
Rogers, L.C.G.: It\^o excursion theory via resolvents. {\it Z. Wahrsch. verw.
Gebiete} {\bf 63} (1983) 237--255.
\smallskip

\itemitem{[\SAL]}
Salminen, P.: One-dimensional diffusions and their exit spaces. {\it Math.
Scand.\/} {\bf 54} (1984) 209--220. 
\smallskip
 
\itemitem{[\SHP]}
Shepp, L.A.: The joint density of the maximum and its location for a Wiener
process with drift. {\it J. Appl. Probab.} {\bf  16} (1979)  423--427.
\smallskip

\itemitem{[\VER]}
Vervaat, W.:  A relation between Brownian bridge and Brownian excursion. {\it
Ann. Probab. } {\bf 7} (1979)  143--149. 
\smallskip

\itemitem{[\WIL]}
Williams, D.:  Path decomposition and continuity of local time for
one-dimensional diffusions, I. {\it Proc. London Math. Soc.} (3) {\bf 28} 
(1974), 738--768. 

\end